\documentclass[a4paper,11pt]{amsart}
\usepackage{amscd,amsmath,amsthm}
\newtheorem{theorem}{Theorem}[section]
\newtheorem{cor}[theorem]{Corollary}
\newtheorem{lemma}[theorem]{Lemma}
\newtheorem{prop}[theorem]{Proposition}
\theoremstyle{remark}
\newtheorem{remark}[theorem]{Remark}
\theoremstyle{remark}
\newtheorem{example}[theorem]{Example}
\theoremstyle{definition}

\numberwithin{equation}{section}
\DeclareMathOperator{\Ind}{Ind}
\DeclareMathOperator{\dashind}{-Ind}
\DeclareMathOperator{\Aut}{Aut}
\DeclareMathOperator{\Ad}{Ad}
\DeclareMathOperator{\clsp}{\overline{span}}
\DeclareMathOperator{\lsp}{{span}}
\newcommand{\abs}[1]{\lvert#1\rvert}
\newcommand{\norm}[1]{\lVert#1\rVert}
\newcommand{\cstar}{$C^*$\ndash}

\newcommand{\tx}{\Tt_X}
\newcommand{\compact}[2]{\Theta_{#1,#2}}
\newcommand{\ip}[2]{\langle #1, #2 \rangle}
\newcommand{\bh}{B(\Hh)}

\newcommand{\ndash}{\nobreakdash-}
\newcommand{\field}[1]{\mathbb{#1}}
\newcommand{\CC}{\field{C}}
\newcommand{\FF}{\field{F}}
\newcommand{\NN}{\field{N}}

\newcommand{\TT}{\field{T}}
\newcommand{\ZZ}{\field{Z}}

\newcommand{\Hh}{{\mathcal H}}
\newcommand{\Kk}{{\mathcal K}}
\newcommand{\Ll}{{\mathcal L}}
\newcommand{\Mm}{{\mathcal M}}
\newcommand{\Oo}{{\mathcal O}}

\newcommand{\Tt}{{\mathcal T}}

\begin{document}
%
% Top Matter
%
\title[The Toeplitz algebra of a Hilbert bimodule]
{The Toeplitz algebra of a Hilbert bimodule}
\author[Neal J. Fowler]{Neal J. Fowler}
\address{Department of Mathematics  \\
      University of Newcastle\\  NSW  2308\\ AUSTRALIA}
\email{neal@math.newcastle.edu.au, iain@maths.newcastle.edu.au}
\author[Iain Raeburn]{Iain Raeburn}
\thanks{This research was supported by the Australian Research Council.}
\date{June 20, 1998}
\begin{abstract}
Suppose a $C^*$-algebra $A$ acts by adjointable operators on a Hilbert
$A$-module $X$.  Pimsner constructed a $C^*$-algebra $\mathcal O_X$ which
includes, for particular choices of $X$, crossed products of $A$ by $\mathbb Z$,
the Cuntz algebras $\mathcal O_n$, and the Cuntz-Krieger algebras $\mathcal O_B$.
Here we analyse the representations of the corresponding Toeplitz algebra.
One consequence is a uniqueness theorem for the Toeplitz-Cuntz-Krieger
algebras of directed graphs, which includes Cuntz's uniqueness theorem for
$\mathcal O_\infty$. 
\end{abstract}
\maketitle
%
% Document body
%

A Hilbert bimodule $X$ over a $C^*$-algebra $A$ is a right Hilbert $A$-module
with a left action of $A$ by adjointable operators. The motivating example  comes 
from an automorphism $\alpha$ of $A$: take $X_A=A_A$, and define the left action
of $A$ by $a\cdot b:=\alpha(a)b$. In \cite{pimsner}, Pimsner constructed a
$C^*$-algebra $\Oo_X$ from a Hilbert bimodule $X$ in such a way that the $\Oo_X$
corresponding to an automorphism $\alpha$ is the crossed product $A\times_\alpha
\ZZ$. He also produced interesting examples of bimodules which do not arise from
automorphisms or endomorphisms, including bimodules over
finite-dimensional commutative $C^*$-algebras for which the corresponding $\Oo_X$
are the Cuntz-Krieger algebras. The Cuntz algebra $\Oo_n$ is $\Oo_X$
when $_\CC X_\CC$  is a Hilbert space of dimension $n$ and the
left action of $\CC$ is by multiples of the identity.

Here we use methods developed in \cite{lacarae,fowrae} for analysing
semigroup crossed products to study Pimsner's
algebras. These methods seem to apply more directly to Pimsner's
analogue  of the Toeplitz-Cuntz algebras rather than his analogue
$\Oo_X$ of the Cuntz algebras. Nevertheless, our results yield new information
about the Cuntz-Krieger algebras of some infinite graphs, giving a whole class
of these algebras which behave like
$\Oo_\infty$.

The uniqueness theorems for $C^*$-algebras generated by algebraic systems
of isometries say, roughly speaking, that all examples of a given system in which
the isometries are non-unitary generate isomorphic $C^*$-algebras. We can approach
such a theorem by introducing a $C^*$-algebra which is universal for 
 systems of the given type, and then characterising its faithful representations.
Here the systems consist of representations $\psi$ of $X$ and $\pi$ of $A$ on the
same Hilbert space which convert the module actions and the inner product to
operator multiplication; we call these Toeplitz representations of $X$. (The
partial isometries and isometries appearing in more conventional systems are
obtained by applying $\psi$ to the elements of a basis for $X$.) In Section 1, we
discuss these Toeplitz representations, show that there is a universal
$C^*$-algebra
$\tx$ generated by a Toeplitz representation, and prove some general results
relating these representations to the induced representations of Rieffel.

Our first theorem is very much in the spirit of other theorems about
$C^*$-algebras generated by systems of isometries: it gives a condition on a
Toeplitz representation $(\psi,\pi)$ which implies that the
corresponding representation
$\psi\times\pi$ of $\tx$ is faithful (Theorem~\ref{theorem:TX}). In broad terms,
this condition says that the ranges of all the operators $\psi(x)$ should leave
enough room for
$A$ to act faithfully. The proof follows standard lines: we use a canonical gauge
action
$\gamma$ to construct an expectation onto a core $\tx^\gamma$, and then show
that both the core and the expectation are implemented faithfully in the given
Toeplitz representation.

When the bimodule $_\CC X_\CC$ is an infinite-dimensional Hilbert space,
Theorem~\ref{theorem:TX} says that a family $\{S_i:i\in \NN\}$ of isometries on
$\Hh$ with orthogonal ranges generates a faithful representation of $\Oo_\infty$
if the ranges $S_i\Hh$ do not span $\Hh$. However, more is true: Cuntz proved
that
\emph{every} family of isometries with orthogonal ranges generates a faithful
representation of $\Oo_\infty$. Our main theorem is an improvement of
Theorem~\ref{theorem:TX} which gives the full strength of Cuntz's result
(Theorem~\ref{theorem:TX2}): we assume that $X$ has a direct-sum
decomposition $X=\bigoplus_\lambda X^\lambda$, but only ask that $A$ acts
faithfully on $(\bigoplus_{\lambda\in F}\psi(X^\lambda)\Hh)^\perp$ for every
finite subset $F$ of indices. For $_\CC X_\CC$, the decomposition is parametrised
by a basis of $X$, and the hypothesis asks that $\sum_{i=1}^n S_iS_i^*<1$ for all
finite $n$, which is trivially true if there are infinitely many $S_i$. To prove
Theorem~\ref{theorem:TX2}, we use the direct-sum decomposition to go further into
the core; we need  the special case in Theorem~\ref{theorem:TX} to construct the
expectation which does this.

The new applications of our theorem involve the $C^*$-algebras of directed
graphs. For a locally finite graph $E$, the $C^*$-algebra $C^*(E)$ is by
definition universal for Cuntz-Krieger $E$-families: families $\{S_f\}$ of
partial isometries, parametrised by the edge set $E^1$ of the graph, and
satisfying in particular $S_e^*S_e=\sum_{\{f:s(f)=r(e)\}}S_fS_f^*$, where
$r,s:E^1\to E^0$ send edges to their range and source vertices
\cite{kprr,kpr}.  The graph algebra $C^*(E)$ can be realised in a very natural
way as the Cuntz-Pimsner algebra $\Oo_X$ of a bimodule $X$ over the algebra
$A=c_0(E^0)$ (see \cite{pinzari,kpw2} and Example~\ref{graphbimod} below). For
graphs in which vertices can emit infinitely many edges, the Cuntz-Krieger
relations involve infinite sums which do not make sense in a $C^*$-algebra, and
it is not clear how to best define a useful notion of graph $C^*$-algebra. We
show that this problem disappears if all vertices emit infinitely many edges:
all families satisfying $S_e^*S_e\geq\sum_{\{f:s(f)=r(e)\}}S_fS_f^*$ generate
isomorphic $C^*$-algebras (Theorem~\ref{theorem:tck}). If the graph is also
transitive, this $C^*$-algebra is simple (Corollary~\ref{oxsimple}). 

Since Hilbert bimodules are a relatively new field of study, and since they arise in so
many different ways, the precise axioms are not yet standard. Thus different authors have
assumed that
$\phi:A\to\Ll(X)$ is injective, that $A$ acts by compact operators on $X$, that
$A$ acts nondegenerately on $X$, or that
$X$ is full.  We have been careful to avoid such assumptions, and in our final section we
illustrate using the bimodules of graphs why we believe this to be helpful. We also
give a couple of new applications involving other classes of Hilbert bimodules.

\section{Toeplitz representations and the Toeplitz algebra}

By a {\em Hilbert bimodule\/} over a \cstar algebra $A$
we shall mean a right Hilbert $A$-module $X$ together with an action
of $A$ by adjointable operators  on $X$. The left
action gives a homomorphism  of $A$ into the \cstar algebra $\Ll(X)$ of 
adjointable operators, which we denote by $\phi$.

A {\em Toeplitz representation} $(\psi,\pi)$ of a Hilbert bimodule
$X$
 in a
\cstar algebra
$B$ consists of a linear map $\psi:X\to B$
and a homomorphism $\pi:A\to B$ such that
\begin{align}
\psi(x\cdot a) &= \psi(x)\pi(a),\label{eq:rep1} \\
\psi(x)^*\psi(y)&= \pi(\ip xy_A),\ \mbox{ and}\label{eq:rep2} \\
\psi(a\cdot x) &= \pi(a)\psi(x).\label{eq:rep3}
\end{align}
for $x,y\in X$ and $a\in A$.
When  $B=\bh$
for some  Hilbert space $\Hh$, we call $(\psi,\pi)$ a Toeplitz
representation of $X$ on $\Hh$.

\begin{remark}\label{remark:psi} In fact Condition \eqref{eq:rep2}
implies that $\psi$ is linear, as in \cite[p.8]{arv}.
It also implies that $\psi$ is bounded: for $x\in X$ we have
\[
\norm{\psi(x)}^2 = \norm{\psi(x)^*\psi(x)} = \norm{\pi(\ip xx_A)}
\le \norm{\ip xx_A} = \norm x^2.
\]
If $\pi$ is injective, then we have equality throughout, and $\psi$
is isometric.
\end{remark}

While many important examples of Hilbert bimodules are given in
\cite[\S1]{pimsner}, \cite[Example 22]{muhly solel2}
and \cite[\S3]{muhly solel}, the examples of most interest to us are
associated to an infinite directed graph. These are not entirely new:  it is
shown in
\cite[p.193]{pimsner} how to build a bimodule from a finite
$\{0,1\}$-matrix $A$, and that bimodule can be obtained by applying
the following construction to the finite graph with incidence matrix $A$.
However, the simplicity of the formulas in the next
Example suggests that it may be more natural to think in terms of
graphs rather than $\{0,1\}$-matrices. 

\begin{example}[The Cuntz-Krieger bimodule]\label{graphbimod}
Suppose $E=(E^0,E^1,r,s)$ is a directed graph with vertex set $E^0$,
edge set $E^1$, and $r,s:E^1\to E^0$ describing the range and source
of edges. Let
$X = X(E)$ be the vector space of functions
$x:E^1\to \CC$ for which the function
\[
v\in E^0 \mapsto \sum_{\{f\in E^1: r(f) = v\}} \abs{x(f)}^2
\]
belongs to $A := c_0(E^0)$.  Then with the operations
\begin{align*}
(x\cdot a)(f) &:= x(f)a(r(f)) \ \mbox{ for } f\in E^1,\\
\ip xy_A(v) &:= \sum_{\{f\in E^1: r(f) = v\}}
\overline{x(f)}y(f) \ \mbox{ for }v\in E^0, \mbox{ and}\\
(a\cdot x)(f) &:= a(s(f))x(f) \ \mbox{ for } f\in E^1,
\end{align*}
$X$ is a Hilbert bimodule over
$A$.

Both the module $X$ and the algebra $A$ are spanned in an
appropriate sense by point masses $\delta_f,\delta_v$, and
we have
\[
\ip {\delta_e}{\delta_f}_A=
\begin{cases}
  \delta_{r(e)}
    & \text{if $e=f$} \\
  0
    & \text{otherwise;} \\
\end{cases}
\]
the elements $\delta_v$ are a family of mutually orthogonal
projections in the $C^*$-algebra $A$. If $(\psi,\pi)$ is a Toeplitz
representation of this Hilbert bimodule
$X$ on $\Hh$,  then the operators $P_v:=\pi(\delta_v)$ are mutually
orthogonal projections on $\Hh$, and
\eqref{eq:rep2} implies that the operators $S_f:=\psi(\delta_f)$ are
partial isometries with initial projection $P_{r(f)}$ and
mutually orthogonal range projections; \eqref{eq:rep3} implies that
these range projections satisfy
\begin{equation}\label{eq:tck}
\sum_{\{f\in E^1: s(f) = v\}} S_fS_f^* \le P_v\ \ \mbox{ for } v\in
E^0.
\end{equation}
We say that $\{S_f,P_v\}$ is a \emph{Toeplitz-Cuntz-Krieger
family\/} for the graph $E$. Conversely, given any such family on
$\Hh$, we can define a representation $\pi:A\to B(\Hh)$ by
$\pi(a):=\sum_v a(v)P_v$, and a linear map $\psi:C_c(E^1)\to B(\Hh)$
by $\psi(x):=\sum_f x(f)S_f$; routine calculations show that $\psi$
is isometric for the $A$-norm on $C_c(E^1)\subset X$ and hence
extends to a linear map on all of $X$, and that $(\psi,\pi)$ is a
Toeplitz representation of $X$. 
  \end{example}

\begin{prop}\label{prop:universal}
Let $X$ be a Hilbert bimodule over $A$.
Then there is a \cstar algebra $\tx$
and a Toeplitz representation $(i_X,i_A):X\to \tx$ such that

\textup{(a)} for every Toeplitz representation $(\psi,\pi)$ of $X$,
there is a homomorphism $\psi\times\pi$ of $\tx$ such that
$(\psi\times\pi)\circ i_X = \psi$ and $(\psi\times\pi)\circ i_A = \pi$; and

\textup{(b)} $\tx$ is generated as a \cstar algebra by $i_X(X) \cup i_A(A)$.

\noindent The triple $(\tx,i_X,i_A)$ is unique:
if
$(B,i_X',i_A')$ has similar properties, there is an isomorphism
$\theta:\tx\to B$ such that
$\theta\circ i_X = i_X'$ and $\theta\circ i_A = i_A'$. Both maps
$i_X$ and $i_A$ are injective. There is a strongly continuous action
$\gamma:\TT\to \Aut\tx$ such that $\gamma_z(i_A(a))=i_A(a)$ and
$\gamma_z(i_X(x))=zi_X(x)$ for $a\in A, x\in X$.
\end{prop}

We call $\tx$ the {\em Toeplitz algebra\/} of $X$ and $\gamma$ the
\emph{gauge action}. To prove the existence of $\tx$, we need to know
that the bimodule has lots of nontrivial representations. Here the
fundamental example is a modification of Fock space, due essentially
to  Pimsner
\cite{pimsner}.

\begin{example}[The Fock representation]\label{example:Fock}
For $n\ge 1$, the $n$-fold internal tensor product
$X^{\otimes n} := X\otimes_A\dotsm\otimes_A X$ is naturally a right
Hilbert $A$-module, and $A$ acts on the left by
\[
a\cdot(x_1\otimes_A\dotsm\otimes_A x_n) :=
(a\cdot x_1)\otimes_A\dotsm\otimes_A x_n;
\]
if we need a name for the operator we call it
$\phi(a)
\otimes_A 1^{n-1}$, and we continue to write $x$ for a
typical element of $X^{\otimes n}$. For $n=0$, we take $X^{\otimes 0}$
to be the Hilbert module
$A$ with left action
$a\cdot b: = ab$. Then the Hilbert-module direct sum $F(X) :=
\bigoplus_{n=0}^\infty X^{\otimes n}$ carries a diagonal left action
of $A$ in which $a\cdot(x_n):=(a\cdot x_n)$. We can induce a
representation
$\pi_0:A\to B(\Hh)$ to a representation
$F(X)\dashind_A^{\Ll(F(X))}\pi_0$ of $\Ll(F(X))$ on
$F(X)\otimes_A\Hh$, which restricts to a representation
$\pi:=F(X)\dashind_A^{A}\pi_0$ of $A$. 

For each  $x\in X$, we can
define a \emph{creation operator\/}
$T(x)$  on $F(X)$ by
\[
T(x)y=\begin{cases}
 x\cdot y
    & \text{if $y\in X^{\otimes 0}=A$} \\
  x\otimes_A y
    & \text{if $y\in X^{\otimes n}$ for some $n\geq 1$;} \\
\end{cases}
\]
routine calculations show that $T(x)$ is adjointable with
\[
T(x)^*z=\begin{cases}
 0
    & \text{if $z\in X^{\otimes 0}=A$} \\
  \ip {x}{x_1}_A\cdot y
    & \text{if $z=x_1\otimes_A y\in X\otimes_A X^{\otimes n-1}=
X^{\otimes n}$.}
\\
\end{cases}
\]
If we now define $\psi:X\to B(F(X)\otimes_A\Hh)$ by 
\[\psi(x):=
F(X)\dashind_A^{\Ll(F(X))}\pi_0(T(x)),
\]
then $(\psi,\pi)$ is a
Toeplitz representation of $X$, called the \emph{Fock
representation induced from $\pi_0$\/}. Note that, since $A$ acts
faithfully on $X^{\otimes 0}=A$ and the representation
$F(X)\dashind_A^{\Ll(F(X))}\pi_0$ is faithful whenever $\pi_0$ is,
 the representation
$\pi$ is faithful whenever $\pi_0$ is; by Remark~\ref{remark:psi}, so
is $\psi$.
\end{example}

\begin{remark}
If we denote by $\phi_\infty$ the diagonal embedding of $A$ 
in $\Ll(F(X))$, then $(T,\phi_\infty)$ is a Toeplitz representation of
$X$ in the
\cstar algebra
$\Ll(F(X))$. Pimsner's
Toeplitz algebra of $X$ is by definition the \cstar subalgebra of
$\Ll(F(X))$ generated by
$T(X)\cup\phi_\infty(A)$ \cite[Definition 1.1]{pimsner}, which is
precisely the image of
$\tx$ under
$T\times\phi_\infty$. In Corollary~\ref{cor:pimsner}, we will show
that our Toeplitz algebra is isomorphic to his by proving
that
$T\times\phi_\infty$ is faithful.
\end{remark}

\begin{proof}[Proof of Proposition~\ref{prop:universal}]
Say that a (Toeplitz) representation $(\psi,\pi)$ of $X$ on a Hilbert
space
$\Hh$ is {\em nondegenerate\/} (resp. {\em cyclic\/})
if the \cstar algebra $C^*(\psi,\pi)$ generated by
$\psi(X)\cup\pi(A)$ acts nondegenerately (resp. cyclically).
For an arbitrary representation $(\psi,\pi)$ of $X$,
let $P$ be the orthogonal projection onto
the essential subspace $\Kk := \overline{C^*(\psi,\pi)\Hh}$;
then $(P\psi,P\pi)$  is a nondegenerate representation of $X$ on $P\Hh$,
and $((I-P)\psi,(I-P)\pi)$ is the zero representation.
By the usual Zorn's lemma argument, $\Kk$ decomposes as a direct sum
of subspaces on which $C^*(\psi,\pi)$ acts cyclically.
Hence every representation is the direct sum of a zero representation and a collection
of cyclic representations.

Let $S$ be a set of cyclic representations of $X$ such that every
cyclic representation of $X$ is unitarily equivalent to an element of
$S$.  (It can be shown that such a set $S$ exists by fixing a Hilbert
space $\Hh$ of sufficiently large dimension, and considering only
cyclic representations on subspaces of $\Hh$. The set $S$ is
nonempty because the Fock representations must have nonzero cyclic
summands.) Let 
\[
\Hh := \bigoplus_{(\psi,\pi)\in S} \Hh_{\psi,\pi},
\ \ i_X := \bigoplus_{(\psi,\pi)\in S} \psi,\ \mbox{ and }\ 
i_A := \bigoplus_{(\psi,\pi)\in S} \pi
\]
(the direct sum defining $i_X$ makes sense because every $\psi$ is contractive).
Then
$(i_X,i_A)$ is a representation of $X$ in $\tx := C^*(i_X,i_A)$; (b) is
satisfied by definition, and (a) can be routinely verified.  

The uniqueness
follows by a standard argument, and the maps $i_X$ and
$i_A$ are injective because the Fock representations factor through
$(i_X,i_A)$ by (a). To establish the existence of the gauge
automorphism $\gamma_z$, just note that $(\tx,z i_X,i_A)$ is also
universal, and invoke the uniqueness. The continuity of the gauge
action follows from a straightforward $\epsilon/3$-argument. 
\end{proof}

Whenever a $C^*$-algebra $C$ acts by adjointable operators on a Hilbert
$A$-module, one can use the module to induce representations of $A$ to
representations of
$C$. If the representation $\pi$ of $A$ is half of a Toeplitz representation, we
can realise the induced representation on the Hilbert space of $\pi$:

\begin{prop}\label{prop:rho}
Let $X$ be a right Hilbert $A$-module,
and suppose $(\psi,\pi)$ is a representation of $X$ on $\Hh$;
that is, $\psi:X\to\bh$ is linear, $\pi:A\to B(\Hh)$ is a representation,
and \eqref{eq:rep1} and \eqref{eq:rep2} hold.

\textup{(1)}
There is a unique representation $\rho=\rho^{\psi,\pi}$ of $\Ll(X)$ on
$\Hh$ with essential subspace $\overline{\psi(X)\Hh}:=\clsp\{\psi(x)h:
x\in X, h\in\Hh\}$ such that
\[
\rho^{\psi,\pi}(S)(\psi(x)h) = \psi(Sx)h\ \mbox{ for } S\in\Ll(X),\ x\in
X\mbox{ and }\ h\in\Hh,
\]
and we then have $\rho(\compact xy) = \psi(x)\psi(y)^*$.

\textup{(2)} If $\Kk$ is a subspace of $\Hh$ which is invariant for
$\pi$, then the subspace $\Mm=\overline{\psi(X)\Kk}$ is invariant for
$\rho$.  If
$\pi|_\Kk$ is faithful, so is $\rho|_\Mm$.
\end{prop}

\begin{proof} (1) 
The map $(x,h)\mapsto \psi(x)h$ is bilinear, and  hence there is  a
linear map $U:
X\odot\Hh\to\Hh$ such that $U(x\otimes h)=\psi(x)h$. Since
\begin{align*}
(U(x\otimes h)\mid U(y\otimes k))
&= (\psi(x)h \mid \psi(y)k) = (h \mid \psi(x)^*\psi(y)k)\\
&= (h \mid\pi(\ip xy_A)k)
= (x\otimes h\mid y\otimes k),
\end{align*}
$U$ extends to an isometry from $X\otimes_A\Hh$ to $\Hh$ such that
$U(x\otimes_A h)=\psi(x)h$.
For $S\in\Ll(X)$ we have
\[
U\Ind\pi(S)U^*(\psi(x)h)
 = U\Ind\pi(S)(x\otimes_A h)
 = U(Sx\otimes_A h)
 = \psi(Sx)h,
\]
so we can define $\rho := \Ad U\circ\Ind\pi$.

 If $x,y,z\in X$, then
\[
\rho(\compact xy)\psi(z)
= \psi(x\cdot\ip yz_A)
= \psi(x)\pi(\ip yz_A)
= \psi(x)\psi(y)^*\psi(z),
\]
so $\rho(\compact xy)$ and $\psi(x)\psi(y)^*$ agree on
$\overline{\psi(X)\Hh}$.
If $k$ is orthogonal to $\psi(X)\Hh$, then $\rho(\compact xy)k = 0$,
so we must show that $\psi(x)\psi(y)^*k = 0$.  
But this follows from
$(\psi(x)\psi(y)^*k \mid h)
= (k \mid \psi(y)\psi(x)^*h) = 0$.

(2) The subspace $\Mm$ is invariant for $\rho$ because $\rho(S)(\psi(x)k) =
\psi(Sx)k$. The restriction of $U$ to $X\otimes_A\Kk$ implements a unitary
equivalence between $\Ind\pi|_{X\otimes_A\Kk}$ and $\rho|_\Mm$;
since the first of these is equivalent to $\Ind(\pi|_\Kk)$,
it is faithful if $\pi|_\Kk$ is, and hence so is $\rho|_\Mm$.
\end{proof}

\begin{remark}\label{pi(1)}
The formula $\rho(\compact xy) = \psi(x)\psi(y)^*$ implies that the
representation
$\rho$ is the canonical extension to $M(\Kk(X))=\Ll(X)$ of the map
Pimsner would call
$\pi^{(1)}$; see \cite[page 202]{pimsner}. (We have avoided the
notation $\pi^{(1)}$ because the map depends on both $\psi$ and $\pi$.)
For a representation $(\psi,\pi)$ of $X$ in a $C^*$-algebra $B$, we
can represent $B$ on a Hilbert space and apply the Proposition to
obtain a homomorphism
$\rho^{\psi,\pi}:\Kk(X)\to B$, but it need not extend canonically to
$\Ll(X)$. 
\end{remark}

\begin{prop}\label{prop:tensor reps}
Let $A$ and $B$ be \cstar algebras,
let $X$ and $Y$ be Hilbert bimodules over $A$, and suppose that $\pi:A\to B$ is a
homomorphism which forms part of Toeplitz representations $(\psi,\pi)$ and
$(\mu,\pi)$ of $X$ and $Y$ in $B$. 

\textup{(1)} There is a linear map $\psi\otimes_A\mu$
of the internal tensor product $X\otimes_A Y$ into $B$ which satisfies
\begin{equation}\label{eq:tensor}
\psi\otimes_A\mu(x\otimes_A y) = \psi(x)\mu(y),\qquad x\in X,\ y\in Y,
\end{equation}
and $(\psi\otimes_A\mu,\pi)$ is a Toeplitz representation of $X\otimes_A Y$.

\textup{(2)} Suppose $B=B(\Hh)$. Denote by $S\mapsto S\otimes_A1$ the canonical
homomorphism of
$\Ll(X)$ into $\Ll(X\otimes_A Y)$ given by the left action of $\Ll(X)$ on $X$,
and let $P_{\psi\otimes_A\mu}$ be the projection of $\Hh$ onto
$\overline{\psi\otimes_A\mu(X\otimes_A Y)(\Hh)}$. Then the representations
$\rho^{\psi,\pi}$ and
$\rho^{\psi\otimes_A\mu,\pi}$ of Proposition~\ref{prop:rho} are related by
\[
\rho^{\psi\otimes_A\mu,\pi}(S\otimes_A1)=\rho^{\psi,\pi}(S)P_{\psi\otimes_A\mu}\
\mbox{ for }S\in\Ll(X).
\]
\end{prop}

\begin{proof} Since $(x,y) \mapsto \psi(x)\mu(y)$ is bilinear,
it induces a linear map $\psi\odot\mu$ on the algebraic
tensor product $X\odot Y$.
For any $x,z\in X$ and $y,w\in Y$ we have
\begin{equation}\label{eq:yxzw}
\begin{split}
\mu(y)^*\psi(x)^*\psi(z)\mu(w)
& = \mu(y)^*\pi(\ip xz_A)\mu(w) \\
& = \mu(y)^*\mu(\ip xz_A\cdot w) \\
& = \pi(\ip y{\ip xz_A\cdot w}_A) \\
& = \pi(\ip{x\otimes_A y}{z\otimes_A w}_A).
\end{split}
\end{equation}
Thus for $v = \sum_i x_i\otimes y_i \in X\odot Y$ we have
\begin{multline*}
\norm{\psi\odot\mu(v)}^2
= \norm{\psi\odot\mu(v)^*\psi\odot\mu(v)}
= \Bigl\lVert \sum_{i,j} \mu(y_i)^*\psi(x_i)^*\psi(x_j)\mu(y_j) \Bigr\rVert \\
= \Bigl\lVert\pi\Big( \sum_{i,j} \ip{x_i\otimes_A y_i}{x_j\otimes_A
y_j}_A\Big)
\Bigr\rVert
\le \Bigl\lVert \sum_{i,j} \ip{x_i\otimes_A y_i}{x_j\otimes_A y_j}_A \Bigr\rVert
= \norm v^2,
\end{multline*}
so $\psi\odot\mu$  induces a contractive
linear map $\psi\otimes_A\mu$ on $X\otimes_A Y$.
Routine calculations on elementary tensors show that
$(\psi\otimes_A\mu,\pi)$ is a Toeplitz representation of
$X\otimes_A Y$.

For part (2), note that the vectors $\psi\otimes_A\mu(x\otimes_A
y)h=\psi(x)\mu(y)h$ span a dense subspace of the essential subspace 
$\overline{\psi\otimes_A\mu(X\otimes_A
Y)\Hh}$ of $\rho^{\psi\otimes_A\mu,\pi}$. Thus the calculation
\begin{align*}
\rho^{\psi\otimes_A\mu,\pi}(S\otimes_A1)(\psi(x)\mu(y)h)
&=\psi\otimes_A\mu\big((S\otimes_A1)(x\otimes_A y)\big)h\\
&=\psi(Sx)\mu(y)h\\
&=\rho^{\psi,\pi}(S)(\psi(x)\mu(y)h)
\end{align*}
implies the result.
\end{proof}

\section{Faithful representations}\label{mainthm}

If $(\psi,\pi)$ is a Toeplitz representation of a Hilbert bimodule
$X$ over $A$ on a Hilbert space $\Hh$, then the subspace
\[
\overline{\psi(X)\Hh}:=\clsp\{\psi(x)h: x\in
X,\ h\in\Hh\}
\]
is invariant for $\pi$: $\pi(a)(\psi(x)h)=\psi(a\cdot x)h$. Thus the
complement $(\psi(X)\Hh)^\perp$ is also invariant for $\pi$. Our
first main theorem says that if $\pi$ is faithful on this
complement, then $\psi\times \pi$ is faithful. 

\begin{theorem}\label{theorem:TX}
Let $X$ be a Hilbert bimodule over a \cstar algebra $A$,
and let $(\psi,\pi)$ be a Toeplitz representation of $X$ on a Hilbert
space
$\Hh$. If $A$ acts faithfully on $(\psi(X)\Hh)^\perp$,
then $\psi\times\pi$ is a faithful representation of
$\tx$. If the homomorphism $\phi:A\to\Ll(X)$ describing the left action of $A$ on
$X$ has range in
$\Kk(X)$ and if $\psi\times\pi$ is faithful, then $A$ acts faithfully on
$(\psi(X)\Hh)^\perp$.
\end{theorem}

Before we prove this theorem we
deduce from it that our Toeplitz
algebra is isomorphic to Pimsner's. This implies in particular that
his algebra is universal for Toeplitz representations
\cite[Theorem~3.4]{pimsner}. 

\begin{cor}\label{cor:pimsner} The Fock representation
$T\times\phi_\infty$ of $\tx$ is faithful.
\end{cor}

\begin{proof}
Let $\pi_0$ be a faithful representation of $A$ on $\Hh$, and
consider
\[
(\psi,\pi):=\Big(\big(F(X)\dashind_A^{\Ll(F(X))}\pi_0\big)\circ T,
\big(F(X)\dashind_A^{\Ll(F(X))}\pi_0\big)\circ\phi_\infty\Big),
\]
which is a Toeplitz representation because $(T,\phi_\infty)$ is. 
For each
$n\ge 0$ and $y\in X^{\otimes n}$, we have $\psi(x)(y\otimes_A
h)=(x\otimes_A y)\otimes_A h$; thus
\[
\overline{\psi(X)(F(X)\otimes_A\Hh)}=\big(\textstyle{\bigoplus_{n=1}^\infty}
X^{\otimes n}\big)\otimes_A
\Hh \cong\bigoplus_{n=1}^\infty(X^{\otimes n}\otimes_A \Hh)
\]
has complement $X^{\otimes 0}\otimes_A \Hh=A\otimes_A\Hh=\Hh$. The
restriction of $\pi$ to this subspace is just
$A\dashind_A^A\pi_0=\pi_0$, which is faithful. Thus
Theorem~\ref{theorem:TX} says that
$\psi\times\pi=\big(F(X)\dashind_A^{\Ll(F(X))}\pi_0\big)
\circ(T\times \phi_\infty)$
is faithful, and hence $T\times\phi_\infty$ is too.
\end{proof}

Averaging over the gauge action gives an  expectation $E$ of
$\tx$ onto the fixed-point algebra $\tx^\gamma$:
\[
E(b) := \int_{\TT} \gamma_w(b)\,dw \ \mbox{ for } b\in\tx.
\]
The map $E$ is a positive linear idempotent
of norm one, and is faithful on positive elements in the sense
that
$E(b^*b) = 0\Longrightarrow b=0$. The main step in the proof of 
Theorem~\ref{theorem:TX} is to show that the expectation $E$ is
spatially implemented: there is a compatible expectation
$E_{\psi,\pi}$ of
$\psi\times\pi(\tx)$ onto $\psi\times\pi(\tx^\gamma)$.

\begin{prop}\label{prop:E and fpa}
Let $(\psi,\pi)$ be a Toeplitz representation of $X$ such that $\pi$ is faithful
on $(\psi(X)\Hh)^\perp$.

\textup{(1)} There is a norm-decreasing map $E_{\psi,\pi}$
on
$\psi\times\pi(\tx)$ such that 
\[
E_{\psi,\pi}\circ (\psi\times\pi) =
(\psi\times\pi)\circ E;
\]

\textup{(2)}  $\psi\times\pi$ is
faithful on the fixed-point algebra
$\tx^\gamma$.
\end{prop}

Before we try to construct $E_{\psi,\pi}$ we need to understand what
$E$ does, and for this we need a description of a dense subalgebra of
$\tx$.

Suppose $(\psi,\pi)$ is a Toeplitz representation of $X$ in a \cstar
algebra
$B$. For $n\ge 1$, Proposition~\ref{prop:tensor reps}
gives us a representation $(\psi^{\otimes n},\pi)$ of the
tensor power $X^{\otimes n}:=X\otimes_A \cdots\otimes_A X$
such that $\psi^{\otimes n}(x_1\otimes_A\cdots \otimes_A
x_n)=\psi(x_1)\cdots\psi(x_n)$. We define
$\psi^{\otimes 0} :=
\pi$. When $m \ge 1$, $X^{\otimes m}\otimes_A X^{\otimes n} =
X^{\otimes (m+n)}$ for every $n\ge 0$, and
$\psi^{\otimes m}\otimes_A \psi^{\otimes n} = \psi^{\otimes (m+n)}$.
There is a slight subtlety for $m = 0$: the natural map $a\otimes_A
x\mapsto a\cdot x$ identifies
$X^{\otimes 0} \otimes_A X^{\otimes n}=A\otimes_A X^{\otimes n}$ 
with the essential submodule $A\cdot X^{\otimes n}$
of $X^{\otimes n}$, and then
$\psi^{\otimes 0}\otimes_A \psi^{\otimes n}$ is the restriction
of $\psi^{\otimes n}$ to this submodule.

\begin{lemma}\label{lemma:monomials} With the above notation, we have
\[
\tx=\clsp\{i_X^{\otimes m}(x)i_X^{\otimes n}(y)^*:m,n\ge 0,
x\in X^{\otimes m},y\in X^{\otimes n}\}.
\]
The expectation $E$ is given by
\[
E(i_X^{\otimes m}(x)i_X^{\otimes n}(y)^*)
 = \begin{cases}i_X^{\otimes m}(x)i_X^{\otimes n}(y)^*&\text{if
$m=n$,}\\
0&\text{if $m\not= n$.}\end{cases}
\]
\end{lemma}

\begin{proof} The algebra $\tx$ is spanned by products of
elements $i_X(x)$, $i_A(a)$ and $i_X(y)^*$; given a word in these
generators, we can usually absorb $i_A(a)$'s into $i_X(x)$'s, and
use $i_X(y)^*i_X(x)=i_A(\ip yx_A)$  to cancel any $i_X(y)^*$
appearing to the left of an $i_X(x)$. (This is
\cite[Lemma~3.1]{pimsner}.)  Since
$\gamma_z(i_X^{\otimes m}(x)i_X^{\otimes n}(y)^*)
 = z^{m-n}i_X^{\otimes m}(x)i_X^{\otimes n}(y)^*$,
the second assertion is easy.
\end{proof}

Lemma~\ref{lemma:monomials} implies that the image $\psi\times\pi(\tx)$ is
spanned by elements $\psi^{\otimes m}(x)\psi^{\otimes n}(y)^*$, and
that $E_{\psi,\pi}$ must satisfy
\begin{equation}\label{spatialE}
E_{\psi,\pi}(\psi^{\otimes m}(x)\psi^{\otimes n}(y)^*)
 = \begin{cases}\psi^{\otimes m}(x)\psi^{\otimes n}(y)^*&\text{if
$m=n$,}\\
0&\text{if $m\not= n$.}\end{cases}
\end{equation}
We shall show that  the formal linear extension $E_{\psi,\pi}$ of the map 
defined by (\ref{spatialE}) is norm-decreasing, and hence extends to a
well-defined norm-decreasing map on $\psi\times\pi(\tx)$. We analyse the norm of
an element
$E_{\psi,\pi}(S)$ by showing that the subspaces
$\overline{\psi^{\otimes n}(X^{\otimes n})\Hh}$ form a decreasing chain of
reducing subspaces, in which the differences are large enough to see operators in
each $\Ll(X^{\otimes n})$ faithfully.

\begin{lemma}\label{movingout}
Suppose that $(\psi,\pi)$ is a Toeplitz representation of $X$ on $\Hh$. For
$n\geq 1$, let
$P_n$ denote the projection of $\Hh$ onto $\overline{\psi^{\otimes n}(X^{\otimes
n})\Hh}$, and let $P_0=1$. Write $\rho_n$ for the representation
$\rho^{\psi^{\otimes n},\pi}$ of $\Ll(X^{\otimes n})$ (so that $\rho_0$ is the
extension of $\pi$ on its essential subspace).

\textup{(1)} We have $P_n\geq P_{n+1}$ for all $n\geq 0$, so
$Q_n:=P_n-P_{n+1}$ is also a projection for $n\geq 0$.

\textup{(2)} For every $n\geq 0$, $k\geq 0$ and
$x,y\in X^{\otimes n}$ we have
\begin{align}
\psi^{\otimes n}(x)P_k&=P_{n+k}\psi^{\otimes n}(x),\ \mbox{
and}\label{eq:commutation}\\ 
P_k\psi^{\otimes n}(x)\psi^{\otimes n}(y)^*&=
\psi^{\otimes n}(x)\psi^{\otimes n}(y)^*P_k.\label{commutation2}
\end{align}

\textup{(3)}  If $\pi$ is faithful on
$(\psi(X)\Hh)^\perp$, then each
$\rho_n$ restricts to a faithful representation of $\Ll(X^{\otimes n})$ on
$Q_n(\Hh)$.
\end{lemma}

\begin{proof} For part (1), observe that the vectors 
$\psi^{\otimes n}(z)\psi(w)h=\psi^{\otimes (n+1)}(z\otimes_A w)h$
span the range of $P_{n+1}$, and are clearly in
the range of $P_n$. 

Equation~(\ref{eq:commutation}) is trivially true for $k=0$ and/or $n=0$.
If $k\ge 1$, $n\geq 1$ and $w\in X^{\otimes k}$, then
\[
\psi^{\otimes n}(x)P_k\psi^{\otimes k}(w)
= \psi^{\otimes n}(x)\psi^{\otimes k}(w)
= P_{n+k}\psi^{\otimes n}(x)\psi^{\otimes k}(w),
\]
so $\psi^{\otimes n}(x)P_k$ and $P_{n+k}\psi^{\otimes n}(x)$ agree on $P_k(\Hh)$.
If $f\in P_k(\Hh)^\perp$, then for any
$z\in X^{\otimes n}$, $w\in X^{\otimes k}$ and $h\in\Hh$ we have
\begin{align*}
(\psi^{\otimes n}(x)f \mid \psi^{\otimes n}(z)\psi^{\otimes k}(w)h)
& = (f \mid \pi(\ip xz_A)\psi^{\otimes k}(w)h) \\
& = (f \mid \psi^{\otimes k}(\ip xz_A\cdot w)h) = 0,
\end{align*}
which implies $P_{n+k}\psi^{\otimes n}(x)f = 0$
because the vectors
$\psi^{\otimes n}(z)\psi^{\otimes k}(w)h=\psi^{\otimes (n+k)}(z\otimes_A w)h$
span the range of $P_{n+k}$. This gives
\eqref{eq:commutation}. Both sides of (\ref{commutation2}) reduce to
$\psi^{\otimes n}(x)\psi^{\otimes n}(y)^*$ when $k< n$; for $k\geq n$, 
\eqref{commutation2} follows from two applications of \eqref{eq:commutation}. 

Part (3) is trivial for $n=0$. For $n\geq 1$, we apply
Proposition~\ref{prop:rho}(2): since
$\pi|_{(I-P_1)\Hh}$ is faithful,
$\rho_n$ is a faithful representation of
$\Ll(X^{\otimes n})$  on $\overline{\psi^{\otimes
n}(X^{\otimes n})(1-P_1)\Hh}$. But this
space is precisely
$Q_n(\Hh)$, because
 \eqref{eq:commutation} implies that  $\psi^{\otimes
n}(x)(1-P_1)=(P_n-P_{n+1})\psi^{\otimes n}(x)$.
\end{proof}

\begin{proof}[Proof of Proposition~\ref{prop:E and fpa}] (1) We have to prove
that for every finite sum 
\[
S:=\sum_{j}
       i_X^{\otimes m_j}(x_j) i_X^{\otimes n_j}(y_j)^*
\]
we have $\|\psi\times\pi(E(S))\|\leq\|\psi\times\pi(S)\|$;
equivalently, we have to prove
\[
\Big\|\sum_{\{j:m_j=n_j\}}
       \psi^{\otimes n_j}(x_j) \psi^{\otimes n_j}(y_j)^*\Big\|\leq
\Big\|\sum_{j}
       \psi^{\otimes m_j}(x_j) \psi^{\otimes n_j}(y_j)^*\Big\|.
\]
We know from \eqref{commutation2} that the projections $Q_k$
commute with every summand in $\psi\times\pi(E(S))$.  If
$m>k$, we have $Q_k\psi^{\otimes m}(x)=Q_kP_m\psi^{\otimes m}(x)=0$, and if
$m\leq k$ and $n\leq k$, \eqref{eq:commutation} gives
\[
Q_k\psi^{\otimes m}(x)\psi^{\otimes n}(y)^*Q_k=\psi^{\otimes
m}(x)Q_{k-m}Q_{k-n}\psi^{\otimes n}(y)^*,
\]
which is $0$ unless $m=n$. Let $K:=\max n_j$. Then
$\rho_K(T\otimes_A1^{K-n})=\rho_n(T)P_K$ by
Proposition~\ref{prop:tensor reps}(2),  so we have
\begin{equation}\label{formforPK}
P_K\big(\psi\times\pi(E(S))\big)
=P_K\rho_K\Big(\sum_{\{j:m_j=n_j\}}\Theta_{x_j,y_j}\otimes_A 1^{K-n_j}\Big);
\end{equation}
because $Q_K\rho_K$ is faithful on $\Ll(X^{\otimes K})$ by the previous lemma, it
follows that
\[
\big\|P_K\big(\psi\times\pi(E(S))\big)\big\|
=\big\|Q_K\big(\psi\times\pi(E(S))\big)\big\|.
\]
Since $Q_0+\cdots +Q_{K-1}+P_K=1$, this gives
\begin{align*}
\|\psi\times\pi(E(S))\|
&=\sup\{\|Q_k\big(\psi\times\pi(E(S))\big)\|:0\leq k\leq K\}\\
&=\sup\{\|Q_k\big(\psi\times\pi(E(S))\big)Q_k\|:0\leq k\leq K\}\\
&=\sup\{\|Q_k\big(\psi\times\pi(S)\big)Q_k\|:0\leq k\leq K\}\\
&\leq\|\psi\times\pi(S)\|.
\end{align*}
Thus
$E_{\psi,\pi}$ extends to a norm-decreasing map on $\psi\times\pi(\tx)$,
giving (1).

Next let $R:=\sum_{j}
       i_X^{\otimes n_j}(x_j) i_X^{\otimes n_j}(y_j)^*$ be a typical finite sum
in the core $\tx^\gamma$; such sums are dense because $E$ is continuous and maps
finite sums to finite sums. For $k<K:=\max n_j$,  Proposition~\ref{prop:tensor
reps}(2) implies that
\[
Q_k(\psi\times\pi(R))
=Q_k\rho_k\Big(\sum_{\{j:n_j\leq k\}}\Theta_{x_j,y_j}\otimes_A 1^{k-n_j}\Big),
\]
and hence
\[
\|Q_k(\psi\times\pi(R))\|\leq 
\Big\|\sum_{\{j:n_j\leq k\}}\Theta_{x_j,y_j}\otimes_A 1^{k-n_j}\Big\|.
\]
There is a similar formula for $\|P_K\big(\psi\times\pi(R)\big)\|$ (see
\eqref{formforPK}), so
\begin{align}
\|\psi\times\pi(R)\|&=\max\big\{\|P_K\big(\psi\times\pi(R)\big)\|;
\|Q_k\big(\psi\times\pi(R)\big)\|:0\leq k< K\big\}\notag\\
&\leq\max\Big\{\Big\|\sum_{\{j:n_j\leq k\}}\Theta_{x_j,y_j}\otimes_A
1^{k-n_j}\Big\|:0\leq k\leq K\Big\}\label{normS}
\end{align}
for every Toeplitz representation $(\psi,\pi)$. Applying this to
a faithful representation shows that \eqref{normS} is an upper bound for
$\|R\|$. 

When $\pi$ is faithful on $(\psi(X)\Hh)^\perp$, the representations $Q_k\rho_k$
and $\rho_K$ are faithful too, so we actually have
\begin{equation}
\|\psi\times\pi(R)\|=\max\Big\{\Big\|\sum_{\{j:n_j\leq
k\}}\Theta_{x_j,y_j}\otimes_A 1^{k-n_j}\Big\|:0\leq k\leq
K\Big\}.\label{isometric}
\end{equation}
In particular, this implies that $\|R\|$ is at least \eqref{normS}; since we
have already seen that $\|R\|$ is at most \eqref{normS}, we deduce that
$\|R\|=\eqref{normS}$, and \eqref{isometric} implies that $\psi\times\pi$ is
isometric on the core.
\end{proof}

\begin{proof}[Proof of Theorem~\ref{theorem:TX}]
Suppose $\pi$ is faithful on $(\psi(X)\Hh)^\perp$ and $S\in \ker\psi\times\pi$.
Then by Proposition~\ref{prop:E and fpa}(1) we have
$\psi\times\pi(E(S^*S))=E_{\psi,\pi}(\psi\times\pi(S^*S))=0$, which by
Proposition~\ref{prop:E and fpa}(2) implies that $E(S^*S)=0$. Because $E$ is
faithful, this forces $S^*S=0$ and
$S=0$.

Now suppose that $\phi(A)\subset\Kk(X)$. Proposition~\ref{prop:rho} gives
a homomorphism $\rho^{i_X,i_A}:\Kk(X)\to\tx$ (see Remark~\ref{pi(1)}), and we
claim that, for any Toeplitz representation $(\psi,\pi)$,
\begin{equation}\label{almostcuntz}
\psi\times\pi\big(i_A(a)-\rho^{i_X,i_A}(\phi(a))\big)=
\pi(a)(1-P_1)=\pi(a)\big|_{(\psi(X)\Hh)^\perp}.
\end{equation}
For any rank-one operator $\Theta_{x,y}$ we have
\[
\psi\times\pi\big(\rho^{i_X,i_A}(\Theta_{x,y})\big)=\psi\times\pi(i_X(x)i_X(y)^*)
=\psi(x)\psi(y)^*=\rho^{\psi,\pi}(\Theta_{x,y}),
\]
and hence $(\psi\times\pi)\circ\rho^{i_X,i_A}=\rho^{\psi,\pi}$ on $\Kk(X)$. On the
other hand, since $\rho^{\psi,\pi}(\phi(a))$ agrees with $\pi(a)$ on
$\psi(X)\Hh$, we have $\rho^{\psi,\pi}(\phi(a))=\pi(a)P_1$. These two
observations give the claim \eqref{almostcuntz}.

Since there are Toeplitz representations $(\psi,\pi)$ in which $\pi$ is
faithful on $(\psi(X)\Hh)^\perp$ (for example, the Fock representation induced
from a faithful representation of $A$) and $\psi\times\pi$ is then faithful,
\eqref{almostcuntz} implies that
$\alpha:a\mapsto i_A(a)-\rho^{i_X,i_A}(\phi(a))$ is an injective homomorphism of
$A$ into $\tx$. (Warning: it is crucial here that $\phi(A)\subset\Kk(X)$.) Thus if
$\psi\times\pi$ is faithful, so is the composition with $\alpha$, and
\eqref{almostcuntz} gives the result.
\end{proof}  

\section{Direct sums of Hilbert bimodules}

If $\{X^\lambda:\lambda\in\Lambda\}$ is a family of Hilbert bimodules over the
same $C^*$-algebra $A$, then the algebraic direct sum $X_0$ is a pre-Hilbert
$A$-module with $(x_\lambda)\cdot a:=(x_\lambda\cdot a)$ and $\ip
{(x_\lambda)}{(y_\lambda)}_A :=
\sum_\lambda \ip{x_\lambda}{y_\lambda}_A$. We can therefore complete $X_0$ to
obtain a Hilbert $A$-module $X$, which we denote by
$\bigoplus_{\lambda\in\Lambda} X^\lambda$ (see \cite[Lemma 2.16]{rw}). 
 There is a left action of
$A$ on $X_0$ defined by $a\cdot(x_\lambda):=(a\cdot x_\lambda)$, which we claim
extends to an action of $A$ by adjointable operators on
$\bigoplus X^\lambda$. To see this, note that the left action of $A$ on each
$X^\lambda$ satisfies $\ip{a\cdot x_\lambda}{a\cdot
x_\lambda}_A\leq\|a\|^2\ip{x_\lambda}{x_\lambda}_A$, and since the sum of positive
elements is positive, we deduce that
\[
\ip{a\cdot (x_\lambda)}{a\cdot
(x_\lambda)}_A \leq\|a\|^2\Big(\sum_\lambda\ip{x_\lambda}{x_\lambda}_A\Big)
=\|a\|^2\ip
{(x_\lambda)}{(y_\lambda)}_A.
\]
Thus the map $(x_\lambda)\mapsto a\cdot(x_\lambda)$ is bounded for the norm
on $X_0$ induced by $\ip{\cdot}{\cdot}_A$, and extends to a map on all of $X$,
which is adjointable with adjoint $(x_\lambda)\mapsto a^*\cdot(x_\lambda)$, as
claimed. We have now shown that $X=\bigoplus_{\lambda\in\Lambda} X^\lambda$
is itself a Hilbert bimodule over $A$, which we call 
the
\emph{direct sum} of the Hilbert bimodules $X^\lambda$.

\begin{theorem}\label{theorem:TX2}
Let $\{X^\lambda:\lambda\in\Lambda\}$ be
a family of Hilbert bimodules over a $C^*$-algebra $A$, 
let
$X:= \bigoplus_{\lambda\in\Lambda} X^\lambda$, and let $(\psi,\pi)$ be a
Toeplitz representation of $X$ on a Hilbert space
$\Hh$. If $A$ acts faithfully on $\big(\psi(\bigoplus_{\lambda\in
F}X^\lambda)\Hh\big)^\perp$ 
for every finite subset $F$ of $\Lambda$, then  $\psi\times\pi$
is faithful on $\tx$. If  $A$ acts by compact operators on the left
of each 
$X^\lambda$ and if $\psi\times\pi$ is faithful,
then $\pi$ acts faithfully on every $\big(\psi(\bigoplus_{\lambda\in
F}X^\lambda)\Hh\big)^\perp$.
\end{theorem}

The proof of this Theorem exploits a   grading of $\tx$ by
the free group $\FF_\Lambda$ on $\Lambda$: picking off the $e$-graded piece 
gives an expectation $E^\Lambda$ which goes further into the core
$\tx^\gamma$ than the expectation $E$ used in Section~\ref{mainthm}.  Such 
gradings are usually formalised in terms of a coaction of
$\FF_\Lambda$ on
$\tx$, but because $\FF_\Lambda$ is not amenable, it would  not be 
obvious from such a formalisation that the associated expectation
$E^\Lambda$ is faithful (see, for example, \cite[\S4]{lacarae}). Here we
shall construct the expectation directly using the Fock representation of
$\tx$, which we know is faithful by Corollary~\ref{cor:pimsner}.

First we need some notation. Let $\FF_\Lambda^+$ be the subsemigroup of
$\FF_\Lambda$ generated by $\Lambda$ and the identity $e$. For $s,t\in
\FF_\Lambda^+$, we write $s\leq t$ if $t$ has the form $sr$ for some
$r\in\FF_\Lambda^+$, and we define
\[
s\vee t:=\begin{cases}
t&\text{if $s\leq t$,}\\
s&\text{if $t\leq s$, and}\\
\infty&\text{otherwise.}
\end{cases}
\]
(The pair $(\FF_\Lambda,\FF_\Lambda^+)$ is an example of a {\em
quasi-lattice ordered group\/} (\cite{nica,lacarae}): the
subsemigroup defines a left-invariant partial order on $\FF_\Lambda$ in
which
$s
\le t$ if and only if $s^{-1}t \in
\FF_\Lambda^+$, and, loosely speaking,   every
finite bounded subset has a least
upper bound.)

For a reduced word $s = \lambda_1\dotsm\lambda_n$ in
$\FF_\Lambda^+\setminus\{e\}$, we write $\abs s := n$. We can identify
the tensor power
$X^s := X^{\lambda_1} \otimes_A \dotsm \otimes_A X^{\lambda_n}$
with a submodule of $X^{\otimes n}$. If $(\psi,\pi)$ is a
Toeplitz representation of $X$, we can define
$\psi^\lambda:=\psi|_{X^\lambda}$ and $\psi^s:=\psi^{\lambda_1} \otimes_A
\dotsm \otimes_A \psi^{\lambda_n}$, and then $(\psi^s,\pi)$ is a Toeplitz
representation of $X^s$ by Proposition~\ref{prop:tensor reps}. The associativity
of $\otimes_A$ gives an isomorphism of $X^s\otimes_A X^t$ onto $X^{st}$ which
carries
$\psi^s\otimes_A\psi^t$ into $\psi^{st}$, and that $\psi^s$ agrees with
the restriction of $\psi^{\otimes|s|}$ to $X^s\subseteq X^{\otimes|s|}$.

\begin{prop}\label{prop:covariance}
Let $(\psi,\pi)$ be a Toeplitz representation of $X$ in a $C^*$-algebra.

\textup{(1)} Suppose $s,t\in\FF_\Lambda^+$ and $s\le t$.
Then for every $x,y_1\in X^s$ and $y_2\in X^{s^{-1}t}$ we have
$\psi^s(x)^*\psi^t(y_1\otimes_A y_2) = \psi^{s^{-1}t}(\ip x{y_1}_A\cdot
y_2)$.

\textup{(2)} Suppose $s,t\in\FF_\Lambda^+$ and $s\vee t = \infty$.
Then for every $x\in X^s$ and $y\in X^t$ we have $\psi^s(x)^*\psi^t(y) = 0$.

\textup{(3)} $\psi\times\pi(\tx)=\clsp\{\psi^s(x)\psi^t(y)^*:x\in X^s,
y\in X^t,s,t\in\FF_\Lambda^+\}.$

\textup{(4)} There is a norm-decreasing linear map  $E^\Lambda$  on $\tx$
which satisfies
\[
E^\Lambda(i_X^s(x)i_X^t(y)^*) = \begin{cases}
i_X^s(x)i_X^t(y)^*&\text{if $s=t$ in $\FF_\Lambda^+$,}\\
0&\text{otherwise,}
\end{cases}
\]
and which is faithful on positive elements.
\end{prop}

\begin{proof} Part (1)
is a straightforward computation.
For (2), let $r$ be the longest common initial segment in $s$ and $t$,
so that 
$s = r\lambda s_1$ and
$t = r\mu t_1$ for $r,s_1,t_1 \in\FF_\Lambda^+$ and $\lambda\ne \mu\in
\Lambda$. Then $X^{r\lambda}$ and $X^{r\mu}$ are orthogonal submodules
of $X^{\otimes(|r|+1)}$. Since vectors of the form $x\otimes_A y \in 
X^{r\lambda}\otimes_A X^{s_1}$ span $X^s$ and similarly for $X^t$, the
calculation
\begin{align*}
\psi^s(x\otimes y)^*\psi^t(w\otimes z)
& = \psi^{s_1}(y)^*\psi^{\otimes (|r|+1)}(x)^*\psi^{\otimes
(|r|+1)}(w)\psi^{t_1}(z) \\ & =
\psi^{s_1}(y)^*\pi(\ip{x}{w}_A)\psi^{t_1}(z)
  = 0
\end{align*}
implies (2).

For (3), we show that $C:=\clsp\{\psi^s(x)\psi^t(y)^*\}$ is a \cstar
subalgebra of $\psi\times\pi(\tx)$ which contains $\psi(X)$ and $\pi(A)$.
It is clearly closed under taking adjoints. To see that it is a
subalgebra, consider 
$\psi^s(x)\psi^t(y)^*$ and
$\psi^u(z)\psi^v(w)^*$. Part (2) implies that 
$\psi^t(y)^*\psi^u(z)=0$ if $t\vee u = \infty$. Otherwise,  (1)
implies that  $\psi^t(y)^*\psi^u(z)$ has the form
$\psi^{t^{-1}u}(z')$ (if $t \le u$) or $\psi^{u^{-1}t}(y')^*$ (if $u \le
t$). Absorbing this element into either
$\psi^s(x)$ or $\psi^v(w)^*$ shows that  the product 
$\psi^s(x)\psi^t(y)^*\psi^u(z)\psi^v(w)^*$ belongs to $C$.

Since $X$ is essential as a right $A$-module,
every element has the form $y\cdot a$ for $y\in X$ and $a\in A$.
Approximating $y$ by a finite sum of the form $\sum
y_\lambda$ shows that
$\psi(y\cdot a) = \psi(y)\pi(a) \sim \sum_\lambda
\psi^\lambda(y_\lambda)\psi^e(a^*)^*$
belongs to $C$.
Similarly,  writing an arbitary element of $A$ as $bc^*$ shows that
$\pi(bc^*) = \psi^e(b)\psi^e(c)^*\in C$.

(4) Part (2) implies that the subspaces $X^s$ of $X^{\otimes n}$
corresponding to different words of length $n$ are orthogonal; thus the natural
map is an isomorphism of the Fock bimodule $\bigoplus_{s\in\FF_\Lambda^+} X^s$
onto $F(X)$. For $r\in \FF_\Lambda^+$, let
$R_r$ be the orthogonal projection of $F(X)$ onto $X^r$. Then for each
$S\in\Ll(F(X))$, the sum $\sum_{r\in\FF_\Lambda^+} R_rSR_r$ converges
$*$-strongly to an adjointable operator $\Phi(S)$; the resulting
linear mapping $\Phi$ on $\Ll(F(X))$ 
is idempotent, norm-decreasing, and faithful on
positive operators. Let $T\times\phi_\infty$ be the Fock representation of
$\tx$, which is faithful by Corollary~\ref{cor:pimsner}. We
want to define $E^\Lambda := (T\times\phi_\infty)^{-1} \circ \Phi \circ
(T\times\phi_\infty)$; before we can do this, we need to know that $\Phi$
leaves the range of
$T\times\phi_\infty$ invariant. Both this and the formula in (4) will
follow if we can show that
\begin{equation}\label{expectnicely}
\Phi(T^s(x)T^t(y)^*) = \begin{cases}
T^s(x)T^t(y)^*&\text{if $s=t$ in $\FF_\Lambda^+$,}\\
0&\text{otherwise.}
\end{cases}
\end{equation}

Let
$x\in X^s$ and
$y\in X^t$, and note that $X^r$ is spanned by vectors of the form
$T^r(z)a$, where 
$a\in A = X^{\otimes 0}$.
If $t\vee r = \infty$, then (2)
gives $T^s(x)T^t(y)^*T^r(z)a = 0$.
If $r < t$, then (1) implies that $T^t(y)^*T^r(z)a=T^{r^{-1}t}(y')^*a$ for
some $y'\in X^{r^{-1}t}$, and this vanishes because $T^{r^{-1}t}(y')^*$ kills
$A=X^{\otimes 0}\subset F(X)$. If
$t\le r$, then  
$T^t(y)^*T^r(z)=T^{t^{-1}r}(z')$ for some $z'\in X^{t^{-1}r}$, and
$T^s(x)T^t(y)^*T^r(z)a = T^{st^{-1}r}(x\otimes z')a \in X^{st^{-1}r}$.
Thus
\[
R_rT^s(x)T^t(y)^*R_r =
\begin{cases}
T^s(x)T^t(y)^*R_r&\text{if $s=t$ in $\FF_\Lambda^+$,}\\
0&\text{otherwise,}
\end{cases} 
\]
and summing over $r\in\FF_\Lambda^+$ gives \eqref{expectnicely}.
\end{proof}

Now suppose that $(\psi,\pi)$ is a Toeplitz representation of $X$ on
$\Hh$. As in the previous section, we aim to show that if $\pi$
satisfies 
the hypothesis of Theorem~\ref{theorem:TX2}, then the expectation
$E^\Lambda$ is spatially implemented. The analogues of the projections
$P_n$ are the projections $P_s$ onto the
subspaces $\overline{\psi^s(X^s)\Hh}$, and in the next Lemma we
write down some of their properties. The analogues of the projections
$Q_n=P_n - P_{n+1}$ are the projections $Q_s^F$ described in
Lemma~\ref{partof1}, which is based on \cite[Lemma~1.4]{lacarae};
in Lemma~\ref{lemma:faithful2} we show that $Q_s^F$ is large enough
to see $\Ll(X^s)$ faithfully.

\begin{lemma}\label{propertiesofP}
Let $(\psi,\pi)$ be a Toeplitz representation of $X$ on $\Hh$. For
$s\in
\FF_\Lambda^+$, denote by $\rho_s$ the representation
$\rho^{\psi^s,\pi}:\Ll(X^s)\to B(\Hh)$, and let $P_s$ be the
projection
of $\Hh$ onto $\overline{\psi^s(X^s)\Hh}$; take $P_e=1$ and
$P_\infty=0$.

\textup{(1)} We have $P_sP_t=P_{s\vee t}$ for $s,t\in\FF_\Lambda^+$.

\textup{(2)} For $s,t\in \FF_\Lambda^+$ and $x,y\in X^s$, we have
\begin{align}
\psi^s(x)P_t&=P_{st}\psi^s(x),\ \mbox{
and}\label{commutation3}\\ 
P_t\psi^s(x)\psi^s(y)^*&=
\psi^s(x)\psi^s(y)^*P_t.\label{commutation4}
\end{align}
\end{lemma}

The proofs are like
those of Part~(2) of Lemma~\ref{movingout}; the orthogonality
of $P_s$ and $P_t$ when $s\vee t=\infty$ follows from
Proposition~\ref{prop:covariance}(2).

\begin{lemma}\label{partof1}
Let $F$ be a finite subset of $\FF_\Lambda^+$ such that
$e\in F$. For $s\in F$, let
\[
Q^F_s := P_s\Big(\prod_{\{t\in F: s<t\}} (1 - P_t)\Big).
\]
Then
$1 = \sum_{s\in F} Q^F_s$.
\end{lemma}

\begin{proof}
We proceed by induction on $\abs F$.  If $\abs F = 1$, then $F = \{e\}$,
and $Q^F_e = P_e = 1$.  If $\abs F \ge 2$, we remove a maximal element
$c$ from $F$, and apply the inductive hypothesis to $G :=
F\setminus\{c\}$.  There is a unique longest word $b\in G$ 
such that $b<c$. We claim that only the summand $Q^G_b$ in
the decomposition $1 = \sum_{s\in G} Q^G_s$ is changed by adding $c$ to
$G$; in other words, we claim that
$Q^F_s = Q^{G}_s$ for $s\not= b$. Suppose $s\in G\setminus\{b\}$.  Then
$Q^F_s$ and $Q^G_s$ have the same factors except for an extra
$1-P_c$ in $Q^F_s$ when $s<c$. But   $s<c$ implies $s<b$, because $b$ is
the longest word in $G$ with $b<c$, and $P_bP_c=P_c$ by
Lemma~\ref{propertiesofP}(1); thus $1-P_b=(1-P_b)(1-P_c)$ and $Q^F_s =
Q^{G}_s$, as claimed. 

We now have
$\sum_{s\in F}Q^F_s=\sum_{s\in G\setminus\{b\}}Q^G_s+Q^F_b+Q^F_c$, and it 
suffices to show that
$Q^{G}_b = Q^F_b + Q^F_c$.
If  $t\in G$ and $b<t$, the maximality of $b$  implies that $c \vee t =
\infty$, and hence
$P_cP_t = 0$ by Lemma~\ref{propertiesofP}(1). Thus
\begin{align*}
Q^{G}_b
& = P_b(1-P_c)\Big(\prod_{\{t\in G: b<t\}} (1 - P_t)\Big) +
P_bP_c\Big(\prod_{\{t\in G: b<t\}} (1 - P_t)\Big) \\ & = P_b\Big(\prod_{\{t\in F:
b<t\}} (1 - P_t)\Big) + P_bP_c
   \\
& = Q^F_b + P_c = Q^F_b + Q^F_c,
\end{align*}
as required.
\end{proof}

\begin{lemma}\label{lemma:faithful2}
Suppose $(\psi,\pi)$ is a Toeplitz representation such that 
$A$ acts faithfully on
$\big(\psi(\bigoplus_{\lambda\in F}X^\lambda)\Hh\big)^\perp$ 
for every finite subset $F$ of $\Lambda$.

\textup{(1)} Let $G$ be a finite subset of $\FF_\Lambda^+\setminus\{e\}$,
and let $s\in\FF_\Lambda^+$.  Then $\rho_s:=\rho^{\psi^s,\pi}$
is a faithful representation of $\Ll(X^s)$  on 
$P_s\prod_{t\in G}(1-P_{st})\Hh$.

\textup{(2)} If $F$ is a finite subset of $\FF_\Lambda^+$ with $e\in F$,
then for each $s\in F$, $\rho_s$
is a faithful representation of
$\Ll(X^s)$  on $Q^F_s\Hh$.
\end{lemma}

\begin{proof} (1) Each $t\in G$ has a unique decomposition $t = \lambda_tr$
with $\lambda_t\in\Lambda$ and $r\in\FF_\Lambda^+$;
write $G' := \{\lambda_t:t\in G\}$.
Lemma~\ref{propertiesofP}(1) implies that the
projections $P_\lambda$ for $\lambda\in\Lambda$ are mutually orthogonal,
so $\overline{\psi(\bigoplus_{\lambda\in
G'}X^\lambda)\Hh}=\bigoplus_{\lambda\in G'}P_\lambda\Hh$, and our
hypothesis says that $\pi$ is faithful on the range of
$1 - \sum_{\lambda\in G'} P_\lambda = \prod_{\lambda\in G'} (1 - P_\lambda)$.
But  $P_t \le P_{\lambda_t}$
for each $t$, so
$\prod_{t\in G} (1 - P_t) \ge \prod_{\lambda\in G'} (1 -
P_\lambda)$, and $\pi$ is also faithful on $\prod_{t\in G} (1 - P_t)\Hh$.
Now  Proposition~\ref{prop:rho}(2) implies that $\rho_s$ is faithful on
\[
\Mm_s := \clsp\{\psi^s(x)\big(\textstyle{\prod_{t\in G}} (1 - P_t)h\big): x\in
X^s,\ h\in\Hh\},
\]
which by \eqref{commutation3} is precisely $\prod_{t\in
G}(1-P_{st})P_s\Hh$, at least for $s\ne e$. When $s=e$, $\Mm_e$
is a subspace of $\prod_{t\in G} (1 - P_t)\Hh$, and the result follows.

(2) Apply (1) with $G:=\{s^{-1}t: t\in F,\ s<t\}$.
\end{proof}

We can now construct our spatial implementation of the expectation
$E^\Lambda$.

\begin{prop}\label{prop:E and fpa2}
Suppose $(\psi,\pi)$ is a Toeplitz representation of $\bigoplus X^\lambda$ such
that $A$ acts faithfully on $\big(\psi(\bigoplus_{\lambda\in
F}X^\lambda)\Hh\big)^\perp$ for every finite subset $F$ of $\Lambda$.

\textup{(1)} There is a norm-decreasing linear map
$E^\Lambda_{\psi,\pi}$ on
$\psi\times\pi(\tx)$ such that 
\[
E^\Lambda_{\psi,\pi}\circ (\psi\times\pi)
= (\psi\times\pi)\circ E^\Lambda;
\]

\textup{(2)} $\psi\times\pi$ is faithful on $E^\Lambda(\tx)$.
\end{prop}

\begin{proof}
(1) We show that for each finite sum $S:=\sum_j i_X^{s_j}(x_j)i_X^{t_j}(y_j)^*$,
we have
\[
\Big\|\sum_{\{j:s_j=t_j\}}\psi^{s_j}(x_j)\psi^{s_j}(y_j)^*\Big\|
\leq \Big\|\sum_{j}\psi^{s_j}(x_j)\psi^{t_j}(y_j)^*\Big\|;
\]
then the map
$E^\Lambda_{\psi,\pi}:
\psi\times\pi(S)\mapsto\psi\times\pi(E^\Lambda(S))$
extends to a well-defined norm-decreasing map on
$\psi\times\pi(\tx)$ with the required
properties. 

Let $F:=\{e\}\cup\{s_j\}\cup\{t_j\}$. Equation~\eqref{commutation4}
implies that the projections $P_s$ and $Q^F_s$ commute with every
summand in $\psi\times\pi(E^\Lambda(S))$;
it follows from Lemma~\ref{partof1} that there exists $c\in F$ such that
\[
\|\psi\times\pi(E^\Lambda(S))\|=
\|Q^F_c(\psi\times\pi(E^\Lambda(S)))\|.
\]
If $t\in F$ with $c<t$, then $Q^F_c\psi^t(x)=Q^F_c(1-P_t)P_t\psi^t(x)=0$, and if
$c\vee t=\infty$, then $Q^F_c\psi^t(x)=Q^F_cP_cP_t\psi^t(x)=0$; thus
compressing by $Q^F_c$ kills all summands in $\psi\times\pi(S)$ except possibly those
for which $s_j\leq c$ and $t_j\leq c$.
As in the proof of
Proposition~\ref{prop:E and fpa}, it follows from
Proposition~\ref{prop:tensor reps}(2) that
\[
Q^F_c(\psi\times\pi(E^\Lambda(S)))=
Q^F_c\rho_c\Big(\sum_{\{j:s_j = t_j \leq c\}}\Theta_{x_j,y_j}\otimes_A
1^{s_j^{-1}c}\Big),
\]
and from Lemma~\ref{lemma:faithful2}(2) that 
\[
\|\psi\times\pi(E^\Lambda(S))\|=
\Big\|\sum_{\{j:s_j = t_j \leq c\}}\Theta_{x_j,y_j}\otimes_A
1^{s_j^{-1}c}\Big\|.
\]
The idea now is to replace
$Q^F_c$ by a smaller projection $Q$, in such a way that compressing by
$Q$ kills the remaining off-diagonal terms of $Q^F_c(\psi\times\pi(S))Q^F_c$
but still preserves the norm
of $\psi\times\pi(E^\Lambda(S))$.

For each $s,t\in F$ such that $s \ne t$, $s,t\le c$ and
$s^{-1}c \vee t^{-1}c < \infty$,
we define $d_{s,t}\in \FF_\Lambda^+$ as in \cite[Lemma~3.2]{lacarae}:
\[
d_{s,t} =
\begin{cases}
  (s^{-1}c)^{-1}(t^{-1}c)
    & \text{if $s^{-1}c < t^{-1}c$} \\
  (t^{-1}c)^{-1}(s^{-1}c)
    &  \text{if $t^{-1}c < s^{-1}c$,}\\
\end{cases}
\]
noting in particular that $d_{s,t}$ is never the identity in $\FF_\Lambda^+$.
Let
\[
G := \{c^{-1}t: t\in F,\,c < t\}
     \cup
     \{d_{s,t}\},
\]
and define $Q:= P_c \prod_{t\in G} (I - P_{ct})$.
Notice that we have added factors to the formula for $Q^F_c$, so $Q\leq
Q^F_c$.

To see that $Q$ has the required properties, fix $s,t\in F$ satisfying
$s\not=t$, $s\leq c$ and $t\leq c$. Then from \eqref{commutation3} we
have
\[
Q\psi^s(x)\psi^t(y)^*Q=QP_c\psi^s(x)\psi^t(y)^*P_cQ=
Q\psi^s(x)P_{s^{-1}c}P_{t^{-1}c}\psi^t(y)^*Q,
\]
which certainly vanishes if $s^{-1}c\vee t^{-1}c=\infty$. But if
$s^{-1}c\vee t^{-1}c<\infty$, then $Q\leq P_c-P_{cd_{s,t}}$, so
\begin{align*}
Q\psi^s(x)\psi^t(y)^*Q
&=Q(P_c-P_{cd_{s,t}})\psi^s(x)\psi^t(y)^*(P_c-P_{cd_{s,t}})Q\\
&=Q\psi^s(x)(P_{s^{-1}c}-P_{s^{-1}cd_{s,t}})
(P_{t^{-1}c}-P_{t^{-1}cd_{s,t}})\psi^t(y)^*Q,
\end{align*}
which vanishes because either $s^{-1}cd_{s,t}=t^{-1}c$ or
$t^{-1}cd_{s,t}=s^{-1}c$. We deduce that
\[
Q(\psi\times\pi(S))Q=Q\rho_c
\Big(\sum_{\{j:s_j = t_j \leq c\}}\Theta_{x_j,y_j}\otimes_A
1^{s_j^{-1}c}\Big).
\]
Since $Q\rho_c$ is faithful by Lemma~\ref{lemma:faithful2}(1), we have 
\begin{align*}
\|\psi\times\pi(E^\Lambda(S))\|&=
\Big\|\sum_{\{j:s_j=t_j\leq c\}}\Theta_{x_j,y_j}\otimes_A
1^{s_j^{-1}c}\Big\|\\
&=
\Big\|Q\rho_c\Big(\sum_{\{j:s_j=t_j\leq c\}}\Theta_{x_j,y_j}\otimes_A
1^{s_j^{-1}c}\Big)\Big\|\\
&=\|Q(\psi\times\pi(S))Q\|\\
&\leq\|\psi\times\pi(S)\|,
\end{align*}
giving (1).

Applying the argument of Proposition~\ref{prop:E and fpa}(2) to the
partition $\{Q^F_s\}$ of $1$ gives (2).
\end{proof}

\begin{proof}[Proof of Theorem~\ref{theorem:TX2}]
The first part follows from Proposition~\ref{prop:E and fpa2} just as
Theorem~\ref{theorem:TX} follows from Proposition~\ref{prop:E and fpa}.
Suppose $A$ acts by compact operators on each summand $X^\lambda$. Then $A$
acts by compact operators on $\bigoplus_{\lambda\in F}X^\lambda$ for any
finite set $F$ of indices, giving maps $\phi_F:A\to
\Kk(X)$. An argument like that in the proof of
Theorem~\ref{theorem:TX} shows that
\[
\psi\times\pi\big(i_A(a)-\rho^{i_X,i_A}(\phi_F(a))\big)
=\pi(a)|_{(1-\sum_{\lambda\in F}P_\lambda)\Hh}.
\]
Applying this with $(\psi,\pi)$ satisfying the hypothesis of the
first part implies that
$\alpha_F:a\mapsto i_A(a)-\rho^{i_X,i_A}(\phi_F(a))$ is an injection of
$A$ in $\tx$. If now $(\psi,\pi)$ is a Toeplitz representation for which
$\psi\times\pi$ is faithful, then composing with $\alpha_F$ shows that
the hypothesis is necessary.
\end{proof}

\section{The Toeplitz algebra of a directed graph}\label{section:graphs}

Let $E = (E^0,E^1,r,s)$ be a directed graph and $X(E)$ the Hilbert bimodule over
$A=c_0(E^0)$ discussed in Example~\ref{graphbimod}. Recall that $X(E)$ consists of
functions on the edge set $E^1$, and that $X(E)$ and $A$ are spanned by point
masses $\{\delta_f:f\in E^1\}$ and $\{\delta_v:v\in E^0\}$, respectively.

\begin{theorem}\label{theorem:tck}
The Toeplitz algebra $\Tt_{X(E)}$ is generated by a Toeplitz-Cuntz-Krieger
$E$-family
$\{i_X(\delta_f),i_A(\delta_v):f\in E^1,v\in E^0\}$. It is universal for such
families: if
$\{S_f,P_v\}$ is a Toeplitz-Cuntz-Krieger
$E$-family on a Hilbert space $\Hh$, there is a representation
$\pi^{S,P}:\Tt_{X(E)}\to B(\Hh)$ such that $\pi^{S,P}(i_X(\delta_f))=S_f$ and
$\pi^{S,P}(i_A(\delta_v))=P_v$. The representation $\pi^{S,P}$ is faithful if and
only if every $P_v$ is nonzero (and hence every $S_f$ is nonzero), and
\[
P_v>
\textstyle{\sum_{\{f\in E^1: s(f) = v\}}} S_fS_f^*
\]
for every vertex $v$ which emits at most finitely many edges.
\end{theorem}

\begin{proof}
Write $X := X(E)$.
We proved in Example~\ref{graphbimod} that $\{\psi(\delta_f),\pi(\delta_v)\}$ is
a Toeplitz-Cuntz-Krieger
$E$-family for any Toeplitz representation $(\psi,\pi)$, and this applies in
particular to the canonical representation $(i_X,i_A)$ in $\tx$. The
family generates $\tx$ because $i_X(X)$ and $i_A(A)$ do, because
$\delta_f$ and $\delta_v$ span dense subspaces of $X$ and $A$, and because
$i_X$ and $i_A$ are isometric. We saw in Example~\ref{graphbimod} how the
family $\{S_f,P_v\}$ generates a Toeplitz representation $(\psi,\pi)$ with
$\psi(\delta_f)=S_f$ and $\pi(\delta_v)=P_v$, so $\pi^{S,P}:=\psi\times\pi$ has
the required property. 

For the final statement, we apply Theorem~\ref{theorem:TX2}. For each $f\in
E^1$, we let $X_f$ be the bimodule $\CC$ in which $a\cdot z=a(s(f))z$, $z\cdot
a=za(r(f))$ and $\ip zw_A=\bar zw\delta_{r(f)}$, and note that $(z_f)\mapsto
\sum_f z_f\delta_f$ induces an isomorphism of $\bigoplus_{f\in E^1}X_f$ onto
$X$. (It is easy to check on the algebraic direct sum that the map is a
bimodule homomorphism which preserves the inner products.)  Since
$\Kk(X_f)=\Ll(X_f)$ for each $f$, $A$ acts by compact operators on each $X_f$,
and Theorem~\ref{theorem:TX2} says that $\pi^{S,P}$ is faithful if and only if $A$
acts faithfully on each $(\bigoplus_{f\in F}\Hh_f)^\perp$, where
$\Hh_f=\pi^{S,P}(i_X(\delta_f)\Hh)=S_f\Hh$. The action of $A=c_0(E^0)$ on any
space is faithful iff every $\delta_v$ acts nontrivially, so $A$ acts faithfully
on $(\bigoplus_{f\in F}\Hh_f)^\perp$ if and only if
\begin{equation*}
0 \ne P_v(1 - \textstyle{\sum_{f\in F}} S_fS_f^*) = P_v -
\textstyle{\sum_{\{f\in F
:s(f) = v\}}} S_fS_f^*.
\end{equation*}
If each $P_w$ is nonzero and $v$ emits infinitely many edges,
this holds since
$P_v \ge \sum_{\{f\in E^1: s(f) = v\}} S_fS_f^*$,
so the result follows.
\end{proof}

\begin{cor}
Let $E$ be a directed graph, and suppose that $\{S_f,P_v\}$ and $\{T_f,Q_v\}$ are
Toeplitz-Cuntz-Krieger $E$-families such that each $P_v$ and $Q_v$ is nonzero,
and such that
\[
P_v>
\textstyle{\sum_{\{f\in E^1: s(f) = v\}}} S_fS_f^*\ \mbox{ and }\ 
Q_v>
\textstyle{\sum_{\{f\in E^1: s(f) = v\}}} T_fT_f^*
\]
for every  vertex $v$ which emits at most finitely many edges. Then there is an
isomorphism $\theta$ of $C^*(S_f,P_v)$ onto $C^*(T_f,Q_v)$ such that
$\theta(S_f)=T_f$ for all $f\in E^1$ and $\theta(P_v)=Q_v$ for all $v\in E^0$. 
\end{cor}

\begin{proof}
Take $\theta:=\pi^{T,Q}\circ(\pi^{S,P})^{-1}$.
\end{proof}

\begin{cor} Let $E$ be a directed graph with at least one edge.  Then $\Tt_{X(E)}$ is simple
if and only if every vertex
emits infinitely many edges and  every pair of vertices are joined
by a finite path.
\end{cor}

\begin{proof} First we show that the hypotheses imply simplicity.
Suppose $\theta$ is a representation of $\Tt_{X(E)}$
with a nontrivial kernel, and let $S_f := \theta(s_f)$ and $P_v := \theta(p_v)$.
Since each vertex emits infinitely many edges,
Theorem~\ref{theorem:tck} implies that  $P_v = 0$ for some $v$.
If $s(f)=v$, then
$S_f = P_vS_f = 0$, and hence $P_{r(f)} = S_f^*S_f = 0$ as well.
Since every pair of vertices are joined by a finite path,
it follows that $P_w = 0$ for every $w\in E^0$.
But then $S_f = S_fS_f^*S_f = S_fP_{r(f)} = 0$ for every $f\in E^1$,
and $\theta = 0$.

Conversely, suppose $\Tt_{X(E)}$ is simple. We prove that we can reach every vertex
from a given vertex $v$ by considering the ideal $\langle p_v \rangle$
generated by $p_v$, which is all of $\Tt_{X(E)}$ by simplicity. 
As usual,  we write
$s_\mu := s_{f_1}\dotsm s_{f_n}$ for a finite path $\mu = f_1\dotsm f_n$, define $s_w :=
p_w$ for each vertex $w$, and verify that 
$\Tt_{X(E)} = \clsp\{s_\mu s_\nu^*\}$. The ideal
$\langle p_v \rangle$ is
spanned by products of the form
$s_\mu s_\nu^*p_v s_\sigma s_\tau^*$, which satisfy
\[
s_\mu s_\nu^*p_v s_\sigma s_\tau^*=
\begin{cases}
s_\mu s_{\sigma'} s_\tau^*&\text{if $s(\nu)=s(\sigma)=v$ and $\sigma=\nu\sigma'$,}\\
s_\mu s_{\nu'}^* s_\tau^*&\text{if $s(\nu)=s(\sigma)=v$ and $\nu=\sigma\nu'$, and}\\
0&\text{otherwise.}\end{cases}
\]
On the other hand, if $r(\mu)=r(\tau)$ can be reached from $v$, say by $\alpha$, then $s_\mu
s_\tau^*=s_\mu s_\alpha^* s_\alpha s_\tau^*=s_\mu s_\alpha^* p_v s_\alpha s_\tau^*$
belongs to $\langle p_v \rangle$. Thus
\[
\langle p_v \rangle = \clsp\{s_\mu s_\tau^*: r(\mu) = r(\tau) \in H(v)\},
\]
where $H(v)$ is the set of vertices $w$ for which there is
a path from $v$ to $w$.

We want to prove that $H(v)$ is all of $E^0$. Suppose there exists $w\in E^0\setminus H(v)$.
We shall show that $\|p_w-b\|\geq 1$ for all $b\in \langle p_v \rangle$, which contradicts
$\langle p_v \rangle=\Tt_{X(E)}$.
Suppose $b = \sum \lambda_i s_{\mu_i} s_{\tau_i}^*$ is a typical finite
sum in $\langle p_v \rangle$. Let $F$ be the (finite) set of edges which start at $w$ and
are the initial edge of some $\mu_i$. Theorem~\ref{theorem:tck}
implies that the projection $q := p_w - \sum_{f\in F} s_fs_f^*$ is nonzero. But
$p_ws_{\mu_i}=0$ unless $s(\mu_i)=w$, and then $s_fs_f^*s_{\mu_i}=s_{\mu_i}$ for the
one $f$ which starts $\mu_i$. Thus
\[
qb=\textstyle{\sum_i}\lambda_ip_ws_{\mu_i}s_{\tau_i}^*-
\textstyle{\sum_i}\lambda_i\big(\textstyle{\sum_{f\in
F}}s_fs_f^*\big)s_{\mu_i}s_{\tau_i}^*=0,
\]
and $\norm{p_w - b} \ge \norm{q(p_w - b)} = \norm q = 1$, as required.

The transitivity we have just proved implies that each vertex $v$ emits at
least one edge. If $v$ emits only finitely many edges, then
 $q := p_v -
\sum_{\{f: s(f) = v\}} s_fs_f^*$ is nonzero by Theorem~\ref{theorem:tck}. However, one can
easily construct Toeplitz-Cuntz-Krieger $E$-families on Hilbert space such that
$P_v=\sum_{\{f:s(f)=v\}}S_fS_f^*$, and then
$q$ would be in the kernel of the corresponding representation of $\Tt_{X(E)}$. Thus each
vertex must emit infinitely many edges.
\end{proof}

 In passing from the Toeplitz algebra $\tx$ to the Cuntz-Pimsner algebra
$\Oo_X$, an important role is played by the ideal $J:=\phi^{-1}(\Kk(X))$; the
theory simplifies when this ideal is either $\{0\}$ or $A$, and authors
have often imposed hypotheses which force $J=A$. (This is done, for example, in
\cite{muhly solel} and  \cite{kpw}.) For the bimodules of graphs, one can
identify the ideal $J$ explicitly.

\begin{prop}\label{prop:J}
Let $X(E)$ be the Hilbert bimodule of a directed graph $E$, and let $\phi:A\to
\Ll(X(E))$ be the homomorphism describing the left action of $A=c_0(E^0)$. Then
\[
\phi^{-1}\big(\Kk(X(E))\big)=\clsp\{\delta_v:\text{$v$ emits at most finitely many
edges}\}.
\]
\end{prop}

\begin{proof}
Write $X:=X(E)$. Since $\Kk(X)$ is an ideal in $\Ll(X)$, $J:=\phi^{-1}(\Kk(X))$ is
an ideal in $A=c_0(E^0)$, and hence has the form 
\[
\{a\in A:\text{$a(w)=0$ for $w\notin F$}\}=\clsp\{\delta_v:v\in F\}
\]
for some subset $F$ of the discrete space $E^0$. So it suffices to see that
$\phi(\delta_v)$ belongs to $\Kk(X)$ iff $v$ emits finitely many edges. If $v$
emits finitely  many edges, then
$\phi(\delta_v)=\sum_{\{f:s(f)=v\}}\Theta_{\delta_f,\delta_f}$ is compact.

Suppose now that $v$ emits infinitely many edges. Since $\lsp\{\delta_f\}$ is
dense in $X$ and $(x,y)\mapsto \Theta_{x,y}$ is continuous, we can approximate any
compact operator on $X$ by a finite linear combination of the form
$K:=\sum_{e,f\in F}\lambda_{e,f}\Theta_{\delta_e,\delta_f}$. But for
any such combination $K$, we can find an edge $g\notin F$ such that $s(g)=v$, and
then
$\Theta_{\delta_e,\delta_f}(\delta_g)=\delta_e\cdot\ip{\delta_f}{\delta_g}_A=0$
for all $e,f \in F$. Thus
\begin{align*}
\|\phi(\delta_v)-K\|
&=\sup\{\|(\phi(\delta_v)-K)(x)\|:\|x\|_A\leq 1\}\\ 
&\geq \|\phi(\delta_v)(\delta_g)-K(\delta_g)\|\\
&=\|\delta_g-0\|=1, 
\end{align*}
and hence $\phi(\delta_v)$ is not compact.
\end{proof}

\begin{cor}\label{oxsimple}
If $E$ is a directed graph in which every vertex emits infinitely many edges,
then the Cuntz-Pimsner algebra $\Oo_{X(E)}$ coincides with the Toeplitz algebra
$\Tt_{X(E)}$, and is simple if and only if $E$ is transitive. 
\end{cor}

\begin{remark}
Since (at least in the absence of sources and sinks) the Cuntz-Pimsner algebra
$\Oo_{X(E)}$ is generated by a  Cuntz-Krieger family for the edge matrix
$B$ of $E$, one might guess that
$\Oo_{X(E)}$ is isomorphic to the Cuntz-Krieger algebra $\Oo_B$ of \cite{el}, and
that this last Corollary follows from \cite[Theorem 14.1]{el}.
This guess is correct, but the connection is nontrivial;
since it concerns Cuntz-Pimsner algebras rather than  Toeplitz algebras, we
shall present the details elsewhere. We note also that our Toeplitz algebra
$\Tt_{X(E)}$ is not the Toeplitz-Cuntz-Krieger algebra $\Tt\Oo_B$ discussed in
\cite{el}: their relations do not imply that the initial projections $P_v$
are mutually orthogonal.
\end{remark}

\section{Concluding remarks}
To see why we have avoided placing additional hypotheses
on our bimodules, consider the Cuntz-Krieger
bimodules  of  graphs.  We want to allow
graphs with infinite valency, so Proposition~\ref{prop:J} shows that 
$A$ will not always act  by compact operators.
We also want to consider graphs with sinks (vertices which emit
no edges) and sources (vertices which receive no edges).
Since $v\in E^0$ is a sink iff
$\delta_v\in c_0(E^0)$ acts trivially on the left of $X(E)$, 
$\phi:A\to \Ll(X)$ may not be injective; since $v$ is a source iff $\delta_v$ is
not in the ideal $\clsp\{\ip xy_A\}$, $X$
need not be full as a right Hilbert module.

Every Cuntz-Krieger bimodule $X=X(E)$ is essential, in the sense that $\clsp A\cdot X = X$,
because $\delta_{s(f)}\cdot\delta_f = \delta_f$ for every $f\in E^1$.
However, the following non-essential submodules arise in analysing the ideal structure of
$\Tt_{X(E)}$. Suppose
$V\subset E^0$ is hereditary in the sense that $r(f)\in V$ whenever
$s(f)\in V$.  Then $I := c_0(V)$ is an ideal in $c_0(E^0)$
such that $I\cdot X(E) \subseteq X(E)\cdot I$, so $X(E)\cdot I$ is a Hilbert $I$-bimodule.
However, if there is an edge $f$ such that $s(f)\notin V$ and $r(f)\in V$, then $\delta_f\in
X(E)\cdot I$ but $a\cdot \delta_f=0$ for all $a\in I$. 

Because our modules may not be essential, we cannot require that the
representations $\pi$ in our Toeplitz representations $(\psi,\pi)$ are nondegenerate:
in the Fock representation
induced from a nondegenerate representation of $A$,
$\pi$ is nondegenerate if and only if
$X$ is essential.
Moreover, the essential subspace of $\pi$ need not be invariant
under $\psi$, so it is not in general possible to reduce to the
nondegenerate case as one typically does when dealing with representations of a $*$-algebra. 
The following Corollary illustrates an extreme case: when the left action is
trivial,
$\psi$ and $\pi$ have orthogonal ranges.
In general, we believe the correct notion of nondegeneracy
for a Toeplitz representation $(\psi,\pi)$ is that
the $C^*$-algebra generated by $\psi(X)\cup\pi(A)$
acts nondegenerately; see the proof of Proposition~\ref{prop:universal}.

\begin{cor}
Suppose the left action of $A$ on $X$ is trivial.

\textup{(1)} $\psi\times\pi$ is faithful if and only if $\pi$ is faithful.
If $A$ is simple, so is $\tx$.

\textup{(2)} $\tx$ is canonically isomorphic to the algebra
\[L(X):=\Kk(X\oplus A)=\left(\begin{matrix} \Kk(X) & X \\ \widetilde X
& A
\end{matrix}\right);
\]
if $X_A$ is full, $L(X)$ is the linking algebra of the imprimitivity bimodule
$_{\Kk(X)}X_A$ (see
\cite[\S3.2]{rw}).
\end{cor}

\begin{proof} (1) If $\psi\times\pi$ is faithful, so is
$(\psi\times\pi) \circ i_A = \pi$. On the other hand, for $a\in
A$ and
$x\in X$ we have
$\pi(a)\psi(x) = \psi(a\cdot x) = 0$,
so $\pi$ acts trivially on $\overline{\psi(X)\Hh}$.
Thus if $\pi$ is faithful it must be faithful on $(\psi(X)\Hh)^\perp$, and
$\psi\times\pi$ is faithful by the Theorem.

(2) The formulas
$\psi(x) := \bigl(\begin{smallmatrix} 0 & x \\ 0 & 0
\end{smallmatrix}\bigr)$ and $\pi(a) := \bigl(\begin{smallmatrix} 0 &
0 \\ 0 & a\end{smallmatrix}\bigr)$ define a Toeplitz representation
of
$X$ in
$L(X)$ such that $\pi$ is faithful and $\psi(X) \cup
\pi(A)$ generates $L(X)$. Now use (1).
\end{proof}

Our next application is a different extension of Cuntz's result on the simplicity of
$\Oo_\infty$: to recover it, take each $X^\lambda={}_\CC \CC_\CC$.

\begin{cor}
Let $X$ be a Hilbert bimodule over a simple \cstar algebra $A$.
If 
$X = \bigoplus_{\lambda\in\Lambda} X^\lambda$ and
the left action of $A$ is nontrivial on infinitely many
summands, then the Toeplitz algebra $\tx$ is simple.
\end{cor}

\begin{proof} If $\psi\times\pi$ is a nonzero representation of $\tx$ on $\Hh$, then the
simplicity of $A$ implies that $\pi$ and $\psi$ are faithful. Since the summands in $X$ are
mutually orthogonal, this implies that the action of $\pi$ in each
$\big(\psi\big(\bigoplus_{\lambda\in F}X^\lambda\big)\Hh\big)^\perp$ is nonzero and hence
faithful. Thus the result follows from Theorem~\ref{theorem:TX2}. 
\end{proof}

Our final application is motivated by Pimsner's realisation of 
crossed products by endomorphisms as $\Oo_X$
for suitable $X$.
Let $\tau$ denote the forward-shift endomorphism on the $C^*$-algebra $c$
of bounded sequences, and let $X := \tau(1)c$ be the Hilbert bimodule
over $c$ in which $x\cdot a := xa$, $\ip xy_c := x^*y$ and
$a\cdot x := \tau(a)x$.  Since the identity operator on $X$
is compact, Theorem~\ref{theorem:TX} applies, and we
recover a theorem of Conway, Duncan and Paterson
\cite{cdp} (see also \cite[Theorem~1.3]{hr}).  Recall that
an element $v$ in a \cstar algebra is a \emph{power partial isometry}
if $v^n$ is a partial isometry for every $n\ge 1$.

\begin{prop} $\tx$ is unital,
$v := i_X(\tau(1))^*$ is a power partial isometry,
and $\tx = C^*(1,v)$.
The pair $(\tx,v)$ has the following universal property:
if $B$ is a unital \cstar algebra and $V\in B$
is a power partial isometry, there is a unital homomorphism
$\tx\to B$ which maps $v$ to $V$.
\end{prop}

\begin{proof}
$i_c(1)$ is an identity for $\tx$, and the calculation
\[
i_c(\tau(a)) = i_X(\tau(1))^*i_X(\tau(a))
= v i_X(a\cdot\tau(1)) = vi_c(a)v^*
\]
shows that $v^nv^{*n} = i_c(\tau^n(1))$ is a projection.
These projections and the identity generate $i_c(c)$;
this and $i_X(x) = v^*i_c(x)$ show that $\tx = C^*(1,v)$.

Suppose $V\in B$ is a power partial isometry.
Since $V^nV^{*n} \ge V^{n+1}V^{*(n+1)}$, there is a
unital homomorphism $\pi_V:c\to B$ which satisfies
$\pi_V(\tau^n(1)) = V^nV^{*n}$.
Define $\psi_V(x) := V^*\pi_V(x)$.  We claim that
$(\psi_V,\pi_V)$ is a Toeplitz representation.
Conditions \eqref{eq:rep1} and \eqref{eq:rep2}
for a Toeplitz representation are easy.  For \eqref{eq:rep3} notice that
$\pi_V(\tau(a)) = V\pi_V(a)V^*$, and recall from
\cite{hw} that the initial and range projections of the powers of $V$
form a commuting family, so that $V^*V\in\pi(c)'$; thus
\begin{align*}
\psi_V(a\cdot x)
& = \psi_V(\tau(a)x)
= V^*\pi_V(\tau(a))\pi_V(x)
= V^*V\pi_V(a)V^*\pi_V(x) \\
& = \pi_V(a)V^*VV^*\pi_V(x)
= \pi_V(a)V^*\pi_V(x)
= \pi_V(a)\psi_V(x),
\end{align*}
as required.
Since $\psi_V\times\pi_V(i_c(1)) = \pi_V(1) = 1$
and $\psi_V\times\pi_V(v) = \psi_V(\tau(1))^* = \pi_V(\tau(1))V = VV^*V = V$,
$\psi_V\times\pi_V$ is the desired map.
\end{proof}

\begin{cor} Let $J_n$ denote the truncated shift on $\CC^n$
(with $J_1 := 0$).  Then $C^*(1, \oplus J_n)$ is the universal
unital $C^*$-algebra generated by a power partial isometry.
\end{cor}

\begin{proof} If $V := \oplus J_n$, then Theorem~\ref{theorem:TX}
implies that $\psi_V\times\pi_V$ is faithful.
\end{proof}

\end{document}